\documentclass[12pt]{amsart}
\usepackage{amsthm}
\usepackage{url}
\input xy
\xyoption {all}

\title{The universal simplicial bundle is a simplicial group}
 \author[D.M. Roberts]{David Michael Roberts}
  \address[David Michael Roberts]
  {School of Mathematical Sciences\\
  University of Adelaide\\
  Adelaide, SA 5005 \\
  Australia}
  \email{david.roberts@adelaide.edu.au}



\theoremstyle{plain}   
\newtheorem{thm}{Theorem}
\newtheorem*{thm*}{\textbf{Theorem}}
\theoremstyle{definition}
\newtheorem{definition}[thm]{Definition}

\newcommand{\into}{\hookrightarrow} 

\newcommand{\Top}{\textbf{\textrm{Top}}}
\newcommand{\Set}{\textbf{\textrm{Set}}}
\newcommand{\Grp}{\textbf{\textrm{Grp}}}
\newcommand{\Gpd}{\textbf{\textrm{Gpd}}}

\DeclareMathOperator{\pr}{pr}
\DeclareMathOperator{\hocolim}{hocolim}

\begin{document} 
 
\begin{abstract}  
The classical universal bundle functor $W:s\Grp(C) \to sC$ for simplicial groups in a category $C$ with finite products lifts to a monad on $s\Grp(C)$. This result extends to simplicial algebras for any Lawvere theory containing that of groups.
\end{abstract} 
\maketitle
\tableofcontents


\section{Introduction}

The present note is motivated by two observations. The first, by Segal \cite{Segal_68}, is that the total space $EG$ of the universal bundle for a (well-pointed) topological group $G$ can be chosen to be a topological group. The easiest way to see this is to pass through the simplicial construction Segal introduced, which from the group $G$ gives a simplicial topological group. The geometric realisation of this simplicial group is then a topological group (using the product in the category of $k$-spaces). Additionally, the group $G$ is a (closed) subgroup of $EG$ and the quotient is one of the standard constructions of the classifying space of a topological group.

The second observation, appearing in \cite{Roberts-Schreiber_08}, is that given a strict 2-group $G$ there is a natural construction of a universal bundle $INN(G)$ for $G$ which is a group-like object in $2\Gpd$. Again, there is an inclusion of groups $G \into INN(G)$. This was proved for 2-groups in $\Set$, but also works for strict 2-groups internal to a finitely complete category.

Given a growing interest in higher gauge theory, derived geometry and higher topos theory, it is natural to consider a generalisation of these results to $\infty$-groups, at least in the first instance as presented by simplicial groups. Because all the constructions involved are very simple, we can work internal to an arbitrary category $C$ with finite products.

The functor $W:s\Grp \to s\Set$, introduced in \cite{MacLane_54}, plays the role in the simplicial world analogous to that $E:\Grp(\Top) \to \Top$ does for universal bundles for topological groups. One can easily see that the construction of $W$ works for simplicial groups in a category $C$ with finite products. The main result here is that $W$ lifts (up to isomorphism) through the forgetful functor $s\Grp(C) \to sC$. Furthermore, not only is this an endofunctor on $s\Grp(C)$, it is a monad, with the unit of the monad being the subgroup inclusion.
The result that $WG$ is a group, at least under the assumption that $C$ has all finite limits, is proved in
\cite{Roberts-Stevenson} in a more conceptual manner.

We can say even more about the status of $WG$ as a universal bundle in categories other than $s\Set$, by recent joint work of Nikolaus, Schreiber and Stevenson \cite{NSS_12}. If $\mathcal{C} = Sh_\infty(S)$ is the $\infty$-topos of $\infty$-sheaves on a site $S$ with a terminal object, then any $\infty$-group object $\mathcal{G}$ in $\mathcal{C}$ is presented by a simplicial group $G$ in $sSh(C)$, and moreover every principal $\infty$-bundle is presented by a pullback (in $sSh(C)$) of the universal bundle $WG \to \overline{W}G$ described here.

For background on simplicial objects and simplicial groups the reader may consult the classic \cite{May}. We shall describe simplicial objects in $\Set$, using elements, but all constructions here are possible in a category with finite products, if we take the definition as using generalised elements.

A remark is perhaps necessary for the history of this result. The main theorem of this note was proved around the time \cite{Roberts-Schreiber_08} was written, but the original version of the notes languished, being referred to in one or two places, themselves until now unpublished work. Thanks are due to Jim Stasheff for encouraging a broader distribution. Urs Schreiber and Danny Stevenson made useful suggestions on a draft.

\section{Results}

We must first present the definitions of the objects we are considering. To start with, we have the classical universal bundle for a simplicial group $G$.

\begin{definition}
The \emph{universal $G$-bundle $WG$} has as its set of $n$-simplices 
$$(WG)_n = G_n \times \ldots \times G_0$$
and face and degeneracy operators 
\begin{align*}
d_0(g_n,\ldots,g_0) &= (d_0g_ng_{n-1},g_{n-2},\ldots,g_0),\\
d_i(g_n,\ldots,g_0) &= (d_ig_n,\ldots,d_1g_{n-i+1},d_0g_{n-i}g_{n-i-1},g_{n-i-2},\ldots,g_0),\quad i > 0\\
s_i(g_n,\ldots,g_0) &= (s_ig_n,\ldots,s_0g_{n-i},id_{G_{n-i}},g_{n-i-1},\ldots,g_0)
\end{align*}
\end{definition}
The simplicial group $G$ acts (on the left) on $WG$ by multiplication on the first factor, and the quotient $(WG)/G$ is denoted $\overline{W}G$. We will not need a detailed description for the present purposes, we only need to note that this quotient exists even if we consider simplicial objects internal to other cateories $C$ without assuming existence of colimits.

We now define a simplicial group $W_{gr}G$ for any simplicial group $G$.
\begin{definition}
The set of $n$-simplices of $W_{gr}G$ is given by
$$
(W_{gr} G)_n = G_n \times \ldots \times G_0.
$$
The face and degeneracy maps are
 \begin{align*}
d_0(g_n,\ldots,g_0) &= (g_{n-1},g_{n-2},\ldots,g_0),\\
d_i(g_n,\ldots,g_0) &= (d_ig_n,\ldots,d_1g_{n-i+1},g_{n-i-1},\ldots,g_0),\quad i > 0\\
s_i(g_n,\ldots,g_0) &= (s_ig_n,\ldots,s_0g_{n-i},g_{n-i},g_{n-i-1},\ldots,g_0).
\end{align*}
If we let the product on $(W_{gr} G)_n$ be componentwise, these face and degeneracy maps are homomorphisms, because those of $G$ are. $W_{gr} G$ is then a simplicial group. The construction is clearly functorial.
\end{definition}

We state the main result of this note, and then prove it in the following section after some observations.
The last section contains some observations which are more open-ended.

\begin{thm*}
The endofunctor $W_{gr}\colon  s\Grp(C) \to s\Grp(C)$ is a lift, up to isomorphism, of the universal bundle functor $W$ through the forgetful functor $s\Grp(C) \to sC$. Moreover, $W_{gr}$ is a monad.
\end{thm*}

One immediate extension of this result is to any Lawvere theory 
$T$ extending that of groups.\footnote{A Lawvere theory \cite{Lawvere} roughly corresponds to any algebraic structure that can be defined using maps between finite (including nullary) products, for example monoids, groups, abelian groups, rings and so on. A Lawvere theory is said to extend the theory of groups if every algebra for that Lawvere theory has an underlying group.} Given such a Lawvere theory, there is a forgetful functor $T\textrm{-}Alg(C) \to \Grp(C)$  from the category of $T$-algebras in a category $C$ to the category of groups in $C$. Similarly, we can consider simplicial $T$-algebras in $C$ (equivalently, $T$-algebras in $sC$). We then have a composite functor
$$
sT\textrm{-}Alg(C) \to s\Grp(C) \stackrel{W}{\to} sC
$$
which by the theorem lifts to $s\Grp(C)$. The construction of $W_{gr}$ is such this composite lifts (on the nose) to a functor 
$$
W_T\colon sT\textrm{-}Alg(C) \to sT\textrm{-}Alg(C).
$$
Given the monad structure maps for $W_{gr}$, one easily sees that $W_T$ is also a monad.

One simple observation which is worth making, given our first motivating fact (from Segal's \cite{Segal_68}), is that for $C$ some subcategory of $\Top$, $G$ is a \emph{closed} sub-simplicial group of $W_{gr}G$. Similarly for other algebras in such a $C$ for a more general Lawvere theory.

Finally, given a finite-product-preserving homotopy colimit functor $sC \to C$, it is clear that the object $\hocolim W_{gr} G \in C$ is actually a group object (equiv. a $T$-algebra). We thus 
have come full circle and recovered Segal's result described above (using the fact geometric realisation is a homotopy colimit for a well-pointed simplicial topological group).

\section{Proof}

\noindent {\bf Proof of theorem:} 
There is an isomorphism between $(WG)_n$ and (the underlying set of) $(W_{gr} G)_n$, given by
\begin{align*}
\Phi_n\colon (WG)_n & \to  (W_{gr} G)_n\\
        (g_n,\ldots,g_0) & \mapsto  (k_n,k_{n-1},\ldots,k_0),
\end{align*}
where the $k_j$ are defined recursively as
$$
k_n = g_n,\qquad k_{j-1} = d_0k_{j}g_{j-1}\ \ (j < n).
$$
One can see that the maps $\Phi_n$ define a map $\Phi$ of simplicial sets by the use of the standard identities for the boundary and degeneracy maps for $G$. One can check the inverse map is
\begin{align*}
\Phi^{-1}\colon W_{gr} G & \to  WG\\
        (h_n,\ldots,h_0) & \mapsto  (h_n,d_0h_n^{-1}h_{n-1},\ldots,d_0h_1^{-1}h_0).
\end{align*}
Thus $W_{gr}$ is an up-to-isomorphism lift of $W$. Note that since $W_{gr}G$ is isomorphic to $WG$ its underlying simplicial set is contractible.

We use the isomorphism $\Phi$ to see how $G$ includes into $W_{gr} G$:
\begin{align*}
G_n \into               (WG)_n      & \stackrel{\Phi}{\to}     (W_{gr} G)_n\\
g_n \mapsto (g_n,1,\ldots,1) &            \mapsto           (g_n,d_0g_n,d_0^2g_n,\ldots,d_0^ng_n)
\end{align*}
Call this homomorphism $\iota_G$. As an aside, from this we can see how $G_n$ is closed in $(W_{gr} G)_n$ when we are working with a subcategory $C \into \Top$, as it is given by the conjunction of the collection of equations $g_{n-i} = d_0^ig_n$ for $i=0,\ldots,n$.

Since the (left) action of $G$ on $W_{gr} G$ is defined via $\Phi$, it is trivial to see that $\Phi$ is a $G$-equivariant isomorphism between free $G$-spaces. This means that $\overline{W}G \simeq (W_{gr} G)/G$; this quotient therefore exists in all categories $C$ with finite products.

Thus far we have an endofunctor
$$
W_{gr} \colon s\Grp(C) \to s\Grp(C),
$$
and a natural transformation
$$
\iota\colon 1_{s\Grp} \to W_{gr}
$$
whose component at $G$ is given by the inclusion $\iota_G \colon G \into W_{gr} G$. 

We now have to prove that $W$ is a monad. For background, see for example \cite{CWM}, chapter VI. Notice that 
\begin{eqnarray*}
(W_{gr}^2 G)_n &=& (W_{gr} G)_n \times (W_{gr} G)_{n-1} \times \ldots (W_{gr} G)_0\\
&=&(G_n \times \ldots \times G_0) \times (G_{n-1} \times \ldots \times G_0) \times \ldots \times (G_0).
\end{eqnarray*}
If $\pr_1\colon (W_{gr} G)_j \to G_j$ denotes projection on the first factor, define the maps
$$
(\mu_G)_n = \pr_1 \times \ldots \times \pr_1 \colon (W_{gr}^2 G)_n \to G_n \times \ldots \times G_0 = (W_{gr} G)_n,
$$
which clearly assemble into a map of simplicial groups
$$
\mu_G\colon W_{gr}^2 G \to W_{gr} G,
$$
and these form the components of a natural transformation
$$
\mu:W_{gr}^2 \to W_{gr}.
$$

Now to show that $W_{gr}$ is a monad we need to check that the following diagrams commute
$$
\xymatrix{
W_{gr}^3 G\ar[rr]^{W_{gr} (\mu_G)} \ar[dd]_{\mu_{W_{gr}G}} & & W_{gr}^2 G \ar[dd]^{\mu}\\
\\
W_{gr}^2 G \ar[rr]_{\mu} &&W_{gr}G 
}
$$
and
$$
\xymatrix{
W_{gr} G\ar[rr]^{W_{gr}(\iota_G)} \ar[ddrr]_{=} & & W_{gr}^2 G \ar[dd]^{\mu} && \ar[ll]_{\iota_{W_{gr}G}} W_{gr}G \ar[ddll]^{=}\\
\\
&&W_{gr}G 
}
$$
This can be done level-wise and is a fairly easy if tedious exercise in indices. This completes the proof.

\section{Postscript}

We can prove that $W_{gr}G$ (and more generally $W_TA$ for a $T$-algebra $A$) is contractible in
the category of simplicial groups in $C$ (resp. in $sT\textrm{-}Alg(C)$) directly, rather than showing
the underlying simplicial object is contractible via the isomorphism $WG \stackrel{\sim}{\to} W_{gr}G$.
The extra degeneracies
\begin{align*}
s_{-1}\colon (W_{gr}G)_n & \to  (W_{gr} G)_{n+1}\\
        (g_n,\ldots,g_0) & \mapsto  (1,g_n,g_{n-1},\ldots,g_0),
\end{align*}
give rise to a contracting homotopy in $s\Grp(C)$ (resp. $sT\textrm{-}Alg(C)$). This allows us to
consider the monad $W_{mon}$ on the category of simplicial \emph{monoids} in $C$, which is defined
in exactly the same way as above. $W_{mon}$ lands in the subcategory of \emph{contractible} 
simplicial monoids, but without the inversion operation we cannot display the isomorphism with
the usual construction of $W$ on the category of simplicial monoids. As such the interpretation of
the quotient $(W_{mon}M)/M$ is not straightforward, especially in the level of generality of this paper, 
where we do not assume the existence of colimits. Even when the required quotient exists, it does
not have a nice interpretation as does the universal bundle $W_{gr}G \to \overline{W}G$.

Finally, we observe that the construction of $W_{mon}$ extends to a monad $W_{T,m}$ on the
category of simplicial $T$-algebras in a finite-product category $C$ where $T$ is a Lawvere theory 
containing a specified\footnote{%
That we need a \emph{specified} monoid operation is clear by the example of the theory of associative, 
unital $k$-algebras for $k$ a field; such $k$-algebras have two underlying monoids.%
}
monoid operation $m$ (i.e. an inclusion $m\colon Th(monoids) \into T$ of the theory
of monoids). This monad takes a simplicial algebra for such a Lawvere theory and returns a 
contractible simplicial $T$-algebra containing the original simplicial $T$-algebra.

\end{document}